\documentclass{amsart}

\usepackage{mathtools}
\usepackage{amsmath}
\usepackage{amssymb}
\usepackage{amsthm}
\usepackage{latexsym}
\usepackage[square,numbers]{natbib}

\author{Efe Gürel}
\address{TÜBİTAK Natural Sciences High School, Kocaeli, 41400, Turkey}
\email{efegurel54@gmail.com}

\newtheorem{theorem}{Theorem}[section]

\newtheorem{lemma}[theorem]{Lemma}



\title[Novel Proofs For Fundamental Identities]
 {Novel Proofs For Fundamental Identities Of Weierstrass Sigma And Jacobi Theta Functions}

\begin{document}
\begin{abstract}
In this note, we describe a general procedure to prove functional equations involving quasi-periodic functions. We give novel proofs for fundamental identities of Weierstrass sigma and Jacobi theta functions. Our method is based on the argument principle rather than the classical approach relying on Liouville's theorem. 
\end{abstract}
\maketitle
\section{Introduction}
Let $\omega_1,\omega_3$ be two non-zero complex numbers with $\mathfrak{Im}(\omega_3/\omega_1)>0$. The Weierstrass sigma function $\sigma(z)=\sigma(z;\omega_1,\omega_3)$ is defined as 
\begin{align*}
    \sigma(z)=z\prod_{(n,m)\neq(0,0)}&\exp\left( \frac{z}{2n\omega_1+2m\omega_3} +\frac{z^2}{2\left(2n\omega_1+2m\omega_3\right)^2}\right)\\
    &\left( 1-\frac{z}{2n\omega_1+2m\omega_3} \right).
\end{align*}
Define the half-period $\omega_2$ by $\omega_1+\omega_2+\omega_3=0$ and let the half-period values be $\eta_j=\sigma'(\omega_j)/\sigma(\omega_j)$ for $j=1,2,3$. The $\sigma$ function satisfies the quasi-periodicity relations
\begin{align*}
    \sigma(z+2\omega_j)=-e^{2\eta_j(z+\omega_j)}\sigma(z).
\end{align*}
The auxiliary sigma functions $\sigma_j$ are given by
\begin{align*}
    \sigma_j(z)=e^{\eta_j z}\frac{\sigma\left( \omega_j-z \right)}{\sigma\left( \omega_j\right)}=e^{-\eta_j z}\frac{\sigma\left( \omega_j+z \right)}{\sigma\left( \omega_j \right)}
\end{align*}
and they are quasi-periodic analogous to the function $\sigma$. For a table depicting quasi-periodicity of the functions $\sigma_j$, see \cite{Schwarz}. Let $\tau\in\mathbb{H}$ be a given complex number and $q=e^{\pi i \tau}$ be the associated nome. The four Jacobi theta functions are given by
\begin{gather*}
    \vartheta_1(z,\tau)=-i\sum_{n=-\infty}^{\infty}(-1)^n q^{\left( n+\frac{1}{2} \right)^2}e^{(2n+1)i z}\\
	\vartheta_{2}(z,\tau)=\sum_{n=-\infty}^{\infty}q^{\left( n+\frac{1}{2} \right)^2}e^{(2n+1)i z}\\
	\vartheta_{3}(z,\tau)=\sum_{n=-\infty}^{\infty}q^{ n^2}e^{2n i z }\\
	\vartheta_{4}(z,\tau)=\sum_{n=-\infty}^{\infty}(-1)^n q^{ n^2}e^{2n i z}.
\end{gather*}
We adopt the theta nullwerte notation $\vartheta_j=\vartheta_j(0,\tau)$ when $\tau$ is clear from the context. The four theta functions are quasi-periodic functions with respect to $\pi$ and $\pi\tau$. Incrementing $z$ by a half-period $\omega_j$ permutes the order of sigma functions. Similarly, incrementing $z$ by a half-period $\pi/2,\pi\tau/2$ or $\pi/2+\pi\tau/2$ permutes the order of theta functions. In fact, sigma and theta functions are equal up to an elementary factor \cite{Schwarz}. Indeed, we have the transformation formulas
\begin{gather*}
    \sigma(z)=\frac{2\omega_1}{\pi\vartheta_1'}e^{\frac{\eta_1z^2}{2\omega_1}}\vartheta_1\left( \frac{\pi z}{2\omega_1},\frac{\omega_3}{\omega_1} \right)\\
    \sigma_1(z)=\frac{1}{\vartheta_2}e^{\frac{\eta_1z^2}{2\omega_1}}\vartheta_2\left( \frac{\pi z}{2\omega_1},\frac{\omega_3}{\omega_1} \right)\\
    \sigma_2(z)=\frac{1}{\vartheta_4}e^{\frac{\eta_1z^2}{2\omega_1}}\vartheta_4\left( \frac{\pi z}{2\omega_1},\frac{\omega_3}{\omega_1} \right)\\
    \sigma_3(z)=\frac{1}{\vartheta_3}e^{\frac{\eta_1z^2}{2\omega_1}}\vartheta_3\left( \frac{\pi z}{2\omega_1},\frac{\omega_3}{\omega_1} \right).
\end{gather*}
The sigma and theta functions arise in many contexts in pure and applied mathematics. The sigma and theta functions satisfy many important properties and functional equations. The most significant ones being so called Weierstrass and Jacobi's fundamental identities. For a discussion of these functions, their general properties and their fundamental identities
we refer the reader to \cite{Chandra,FundamentaNova,CourseD'analyse,Koornwinder,Lawdens,Smith,MathematischeWerke}. In this short note, we present a novel method of proving such identities based on contour integration.

\section{Main Results}
We first describe our method generally and then give examples such as Weierstrass $3$-term identity, Weierstrass fundamental identity, Jacobi's addition formulas and Jacobi fundamental identities. Suppose we are trying prove a relation between quasi-periodic functions where both of the sides are entire functions of a fixed variable. The classical approach is to show both sides satisfy the same quasi-periodicity relations and then to consider the ratio of both sides. Then it is proved that the zeros of the both sides are identical. This shows that the ratio of both sides is an entire elliptic function and thus constant by Liouville's theorem. This constant is calculated by expanding the function near a suitable point. Many examples of this method can be seen in \cite{WhittakerWatson}.
\newline

Our method is based on contour integration and the argument principle. Let us consider an entire quasi-periodic function $\varphi$ satisfying the relations
\begin{align}\label{GeneralQP}
    \varphi(z+\lambda_1)=e^{\alpha_1z+\beta_1}\varphi(z)\qquad \text{and}\qquad \varphi(z+\lambda_2)=e^{\alpha_2z+\beta_2}\varphi(z)
\end{align}
where $\mathfrak{Im}(\lambda_2/\lambda_1)>0$. We must show that $\varphi$ is identical to $0$. Assume that $\varphi$ is non-constant. Consider a period parallelogram $P_\omega$ with vertices $\omega,\omega+\lambda_1,\omega+\lambda_1+\lambda_2,\omega+\lambda_2$. Since $\varphi$ is entire, by the argument principle we have
\begin{align*}
    N_\varphi=\frac{1}{2\pi i} \oint_{P_\omega}\frac{\varphi'(z)}{\varphi(z)}dz
\end{align*}
where $N_\varphi$ denotes the number of zeros of $\varphi$ in any period parallelogram. Decomposing the contour integral, we obtain
\begin{align*}
    N_\varphi&=\frac{1}{2\pi i} \left( \int_{\omega}^{\omega+\lambda_1}+\int_{\omega+\lambda_1}^{\omega+\lambda_1+\lambda_2}+\int_{\omega+\lambda_1+\lambda_2}^{\omega+\lambda_2}+\int_{\omega+\lambda_2}^{\omega} \right)\frac{\varphi'(z)}{\varphi(z)}dz\\
    &=\frac{1}{2\pi i}\int_{\omega}^{\omega+\lambda_1}\frac{\varphi'(z)}{\varphi(z)}-\frac{\varphi'(z+\lambda_2)}{\varphi(z+\lambda_2)}dz+\frac{1}{2\pi i}\int_{\omega}^{\omega+\lambda_2}\frac{\varphi'(z+\lambda_1)}{\varphi(z+\lambda_1)}-\frac{\varphi'(z)}{\varphi(z)}dz.
\end{align*}
Using the quasi-periodicity relations \eqref{GeneralQP}, we get
\begin{align}\label{GeneralN}
    N_\varphi=\frac{\alpha_1\lambda_2-\alpha_2\lambda_1}{2\pi i}.
\end{align}
To contradict this statement, we now manually exhibit the zeros of $\varphi$ in a period parallelogram. Existence of more than $\frac{\alpha_1\lambda_2-\alpha_2\lambda_1}{2\pi i}$ zeros contradicts the statement that $\varphi$ is non-constant. Therefore $\varphi$ is identical to $0$.
\newline

We first illustrate our method with Weierstrass $3$-term identity. We note that Legendre's theorem states that $\eta_1\omega_3-\eta_3\omega_1=\frac{\pi i}{2}$.

\begin{theorem}
    For every $z,a,b,c\in\mathbb{C}$, the following equation holds
    \begin{align*}
        0&=\sigma(z+a)\sigma(z-a)\sigma(b+c)\sigma(b-c)\\
	&+\sigma(z+b)\sigma(z-b)\sigma(c+a)\sigma(c-a)\\
	&+\sigma(z+c)\sigma(z-c)\sigma(a+b)\sigma(a-b).
    \end{align*}
\end{theorem}
\begin{proof}
    If any of the numbers $a,b,c$ differ by an amount of $2n\omega_1+2m\omega_3$ for $n,m\in\mathbb{Z}$, the equation trivially holds. Therefore assume that $a,b,c$ are distinct modulo $(2\omega_1,2\omega_3)$. Define the function
    \begin{align*}
        \varphi_1(z)&=\sigma(z+a)\sigma(z-a)\sigma(b+c)\sigma(b-c)\\
	&+\sigma(z+b)\sigma(z-b)\sigma(c+a)\sigma(c-a)\\
	&+\sigma(z+c)\sigma(z-c)\sigma(a+b)\sigma(a-b).
    \end{align*}
    We claim that the function $\varphi_1(z)$ is quasi-periodic. Indeed, we have
    \begin{align*}
        \varphi_1(z+2\omega_1)&=\sigma(z+a+2\omega_1)\sigma(z-a+2\omega_1)\sigma(b+c)\sigma(b-c)\\
	&+\sigma(z+b++2\omega_1)\sigma(z-b+2\omega_1)\sigma(c+a)\sigma(c-a)\\
	&+\sigma(z+c+2\omega_1)\sigma(z-c+2\omega_1)\sigma(a+b)\sigma(a-b)\\
        &=\left( -e^{2\eta_1(z+a+\omega_1)}\sigma(z+a) \right)\left( -e^{2\eta_1(z-a+\omega_1)}\sigma(z-a) \right)\sigma(b+c)\sigma(b-c)\\
        &+\left( -e^{2\eta_1(z+b+\omega_1)}\sigma(z+b) \right)\left( -e^{2\eta_1(z-b+\omega_1)}\sigma(z-b) \right)\sigma(c+a)\sigma(c-a)\\
        &+\left( -e^{2\eta_1(z+c+\omega_1)}\sigma(z+c) \right)\left( -e^{2\eta_1(z-c+\omega_1)}\sigma(z-c) \right)\sigma(a+b)\sigma(a-b)\\
        &=e^{4\eta_1z+4\eta_1\omega_1}\varphi_1(z).
    \end{align*}
    Similarly, we have $\varphi_1(z+2\omega_3)=e^{4\eta_3z+4\eta_3\omega_3}\varphi_1(z)$. Assume $\varphi_1$ is non-constant. Thus by equation \eqref{GeneralN} and Legendre's theorem, the number of zeros of $\varphi_1(z)$ in a period parallelogram is
    \begin{align*}
        N=\frac{8\eta_1\omega_3-8\eta_3\omega_1}{2\pi i}=2.
    \end{align*}
    But we evidently have $\varphi_1(a)=\varphi_1(b)=\varphi_1(c)=0$. This contradicts the assumption that $a,b,c$ are distinct modulo $(2\omega_1,2\omega_3)$. This completes the proof.
\end{proof}
We now prove Weierstrass' fundamental identities \cite{Schwarz}. Since the other fundamental identities can be deduced by incrementing variables by suitable half-periods, it suffices to prove only one. Given complex numbers $a,b,c$ and $d$, we define the associated variables $a',a''$ by
\begin{gather*}
    2a'=a+b+c+d,\qquad \qquad \qquad 2a''=a+b+c-d\\
    2b'=a+b-c-d,\qquad \qquad \qquad 2b''=a+b-c+d\\
    2c'=a-b+c-d,\qquad \qquad \qquad 2c''=a-b+c+d\\
    2d'=-a+b+c-d,\qquad \qquad \qquad 2d''=a-b-c-d.
\end{gather*}
Then we have the following theorem.
\begin{theorem}
    For every $a,b,c,d\in\mathbb{C}$, the following equation holds
    \begin{align*}
        0&=\sigma\left( a \right)\sigma\left( b \right)\sigma\left( c \right)\sigma\left( d \right)\\
        &+\sigma\left( a' \right)\sigma\left( b' \right)\sigma\left( c' \right)\sigma\left( d' \right)\\
        &+\sigma\left( a'' \right)\sigma\left( b'' \right)\sigma\left( c'' \right)\sigma\left( d'' \right).
    \end{align*}
\end{theorem}
\begin{proof}
    Let us define the function
    \begin{align*}
        \varphi_2(a)&=\sigma\left( a \right)\sigma\left( b \right)\sigma\left( c \right)\sigma\left( d \right)\\
        &+\sigma\left( a' \right)\sigma\left( b' \right)\sigma\left( c' \right)\sigma\left( d' \right)
        \\&+\sigma\left( a'' \right)\sigma\left( b'' \right)\sigma\left( c'' \right)\sigma\left( d'' \right).
    \end{align*}
    Then we have the quasi periodicity
    \begin{align*}
        \varphi_2(a+4\omega_1)&=\sigma\left( a+4\omega_1 \right)\sigma\left( b \right)\sigma\left( c \right)\sigma\left( d \right)\\
        &+\sigma\left( a'+2\omega_1 \right)\sigma\left( b'+2\omega_1 \right)\sigma\left( c'+2\omega_1\right)\sigma\left( d'-2\omega_1 \right)\\
        &+\sigma\left( a''+2\omega_1 \right)\sigma\left( b''+2\omega_1 \right)\sigma\left( c''+2\omega_1 \right)\sigma\left( d''+2\omega_1 \right)\\
        &=e^{4\eta_1a+8\eta_1\omega_1}\sigma\left( a \right)\sigma\left( b \right)\sigma\left( c \right)\sigma\left( d \right)\\
        &+e^{2\eta_1\left(a'+b'+c'-d'\right)+8\eta_1\omega_1}\sigma\left( a' \right)\sigma\left( b' \right)\sigma\left( c' \right)\sigma\left( d' \right)\\
        &+e^{2\eta_1\left(a''+b''+c''+d''\right)+8\eta_1\omega_1}\sigma\left( a'' \right)\sigma\left( b'' \right)\sigma\left( c'' \right)\sigma\left( d'' \right)\\
        &=e^{4\eta_1a+8\eta_1\omega_1}\varphi_2(a).
    \end{align*}    
    Similarly, we have $\varphi_2(a+4\omega_3)=e^{4\eta_3a+8\eta_3\omega_3}\varphi_2(a)$. Assume $\varphi_2$ is non-constant. Thus by equation \eqref{GeneralN}, the number of zeros of $\varphi_2(a)$ in a period parallelogram is 
    \begin{align*}
        N=\frac{16\eta_1\omega_3-16\eta_3\omega_1}{2\pi i}=4.
    \end{align*}
    With a similar investigation of the roots of $\varphi_2$, we obtain a contradiction. This completes the proof.
\end{proof}
Finally, we prove a mixed-type sigma function identity involving functions $\sigma_j$. Let us denote $\wp(z)=-\left(\log\sigma(z)\right)''$ with the half-period values $e_j=\wp(\omega_j)$ for $j=1,2,3$.

\begin{theorem}
    Let $k,\ell\in\{1,2,3\}$. Then, for every $z\in\mathbb{C}$, the following equation holds
    \begin{align*}
        \sigma_k^2(z)-\sigma_\ell^2(z)+(e_k-e_\ell)\sigma^2(z)=0.
    \end{align*}
\end{theorem}
\begin{proof}
    Let us define the function
    \begin{align*}
        \varphi_3(z)=\sigma_k^2(z)-\sigma_\ell^2(z)+(e_k-e_\ell)\sigma^2(z).
    \end{align*}
    It has been noted in \cite{Schwarz} that for $j=1,2,3$, we have the quasi-periodicity
    \begin{align*}
        \sigma_k^2(z+2\omega_j)=e^{4\eta_j z+4\eta_j \omega_j}\sigma_k^2(z).
    \end{align*}
    The sigma function $\sigma$ also satisfies the identical quasi-periodicity relation. Therefore $\varphi_3(z+2\omega_1)=e^{4\eta_1z+4\eta_1\omega_1}\varphi_3(z)$ and $\varphi_3(z+2\omega_3)=e^{4\eta_3z+4\eta_3\omega_3}\varphi_3(z)$ . Assume $\varphi_3$ is non-constant. Thus by equation \eqref{GeneralN} and Legendre's theorem, the number of zeros of $\varphi_3(z)$ in a period parallelogram is
    \begin{align*}
        N=\frac{8\eta_1\omega_3-8\eta_3\omega_1}{2\pi i}=2.
    \end{align*}
    Since $\sigma_k(0)=\sigma_\ell(0)=1$ and $\sigma(0)=\sigma_j(\omega_j)=0$ for $j=1,2,3$, we have $\varphi_3(0)=0$. We now prove $\varphi_3(\omega_k)=\varphi_3(\omega_\ell)=0$. It is a classical result in the theory of elliptic functions that
    \begin{align*}
        \sqrt{\wp(z)-e_\ell}=\frac{\sigma_\ell(z)}{\sigma(z)}.
    \end{align*}
    Taking $z=\omega_k$, we obtain $\varphi_3(\omega_k)=0$ and similarly $\varphi_3(\omega_\ell)=0$. This contradicts $N=2$ and completes the proof.
\end{proof}

We now turn to identities of Jacobi theta functions. We first need the following fundamental lemmas about quasi-periodicity \cite{WhittakerWatson}. We suppress the variable $\tau$ for convenience.
\begin{lemma}\label{ThetaQPLemma}
    The following quasi-periodicity relations hold
    \begin{gather*}
        \vartheta_1(z+\pi)=-\vartheta_1(z), \qquad\qquad \qquad \vartheta_1(z+\pi\tau)=-q^{-1}e^{-2iz}\vartheta_1(z)\\
        \vartheta_2(z+\pi)=-\vartheta_2(z), \qquad\qquad \qquad \vartheta_2(z+\pi\tau)=q^{-1}e^{-2iz}\vartheta_2(z)\\
        \vartheta_3(z+\pi)=\vartheta_3(z), \qquad\qquad \qquad \vartheta_3(z+\pi\tau)=q^{-1}e^{-2iz}\vartheta_3(z)\\
        \vartheta_4(z+\pi)=\vartheta_4(z), \qquad\qquad \qquad \vartheta_4(z+\pi\tau)=-q^{-1}e^{-2iz}\vartheta_4(z).
    \end{gather*}
\end{lemma}
\begin{lemma}\label{ThetaTransformLemma}
    The following transformation formulas hold
    \begin{gather*}
        \vartheta_1(z)=-\vartheta_2\left( z+\frac{\pi}{2} \right)=-iq^{1/4}e^{iz}\vartheta_3\left( z+\frac{\pi}{2}+\frac{\pi\tau}{2} \right)=-iq^{1/4}e^{iz}\vartheta_4\left( z+\frac{\pi\tau}{2} \right)\\
        \vartheta_2(z)=q^{1/4}e^{iz}\vartheta_3\left( z+\frac{\pi\tau}{2} \right)=q^{1/4}e^{iz}\vartheta_4\left( z+\frac{\pi}{2}+\frac{\pi\tau}{2} \right)=\vartheta_1\left( z+\frac{\pi}{2} \right)\\
        \vartheta_3(z)=\vartheta_4\left( z+\frac{\pi}{2} \right)=q^{1/4}e^{iz}\vartheta_1\left( z+\frac{\pi}{2}+\frac{\pi\tau}{2} \right)=q^{1/4}e^{iz}\vartheta_2\left( z+\frac{\pi\tau}{2} \right)\\
        \vartheta_4(z)=-iq^{1/4}e^{iz}\vartheta_1\left( z+\frac{\pi\tau}{2} \right)=iq^{1/4}e^{iz}\vartheta_2\left( z+\frac{\pi}{2}+\frac{\pi\tau}{2} \right)=\vartheta_3\left( z+\frac{\pi}{2} \right).
    \end{gather*}
\end{lemma}
We now prove two of the classical addition theorems of Jacobi of different kinds. One with a multiplier of form $\vartheta_j^2$ and one with a multiplier of form $\vartheta_j\vartheta_{j'}$. Similar formulas can be proved in an analogous manner or by incrementing the variables by suitable half-periods. We note that the functions $\vartheta_1(z),\vartheta_2(z),\vartheta_3(z),\vartheta_4(z)$ have zeros at points conguent to $z=0,\frac{\pi}{2},\frac{\pi}{2}+\frac{\pi\tau}{2},\frac{\tau}{2}$ modulo $(\pi,\pi\tau)$ respectively.

\begin{theorem}
    The following equation holds for all $a,b\in\mathbb{C}$
    \begin{align*}
        \vartheta_3(a+b)\vartheta_3(a-b)\vartheta_3^2=\vartheta_3^2(a)\vartheta_3^2(b)+\vartheta_1^2(a)\vartheta_1^2(b).
    \end{align*}
\end{theorem}
\begin{proof}
    If $b$ is congruent to $0$ or $\frac{\pi}{2}+\frac{\pi\tau}{2}$ modulo $(\pi,\pi\tau)$, then the theorem holds trivially by Lemma \ref{ThetaTransformLemma}. Therefore assume otherwise. Let us define the function
    \begin{align*}
        \psi_1(a)=\vartheta_3(a+b)\vartheta_3(a-b)\vartheta_3^2-\vartheta_3^2(a)\vartheta_3^2(b)-\vartheta_1^2(a)\vartheta_1^2(b).
    \end{align*}
    By Lemma \ref{ThetaQPLemma}, we have $\psi_1(a+\pi)=\psi_1(a)$ and
    \begin{align*}
        \psi_1(a+\pi\tau)&=\vartheta_3(a+b+\pi\tau)\vartheta_3(a-b+\pi\tau)\vartheta_3^2
        -\vartheta_3^2(a+\pi\tau)\vartheta_3^2(b)
        -\vartheta_1^2(a+\pi\tau)\vartheta_1^2(b)\\
        &=q^{-2}e^{-4ia}\psi_1(a).
    \end{align*}
    Assume that $\psi_1$ is non-constant. Then by equation \eqref{GeneralN}, the number of zeros of $\psi_1(a)$ in a period parallelogram is $N=2$. We have $\psi_1(0)=0$ and furthermore by Lemma \ref{ThetaTransformLemma},
    \begin{align*}
        \psi_1\left( \frac{\pi}{2}+\frac{\pi\tau}{2} \right)&=\vartheta_3\left( b+\frac{\pi}{2}+\frac{\pi\tau}{2} \right)\vartheta_3\left( \frac{\pi}{2}+\frac{\pi\tau}{2}-b \right)\vartheta_3^2-\vartheta_1^2\left( \frac{\pi}{2}+\frac{\pi\tau}{2} \right)\vartheta_1^2(b)\\
        &=0.
    \end{align*}
    Similarly,
    \begin{align*}
        \psi_1\left(b+ \frac{\pi}{2}+\frac{\pi\tau}{2} \right)&=-\vartheta_3^2\left(b+ \frac{\pi}{2}+\frac{\pi\tau}{2} \right)\vartheta_3^2(b)-\vartheta_1^2\left(b+ \frac{\pi}{2}+\frac{\pi\tau}{2} \right)\vartheta_1^2(b)\\
        &=0.
    \end{align*}
    But since $\psi_1$ only has two zeros in a period parallelogram, we must have that $b$ is congruent to either $0$ or $\frac{\pi}{2}+\frac{\pi\tau}{2}$ modulo $(\pi,\pi\tau)$. Thus contradicting our assumptions. This completes the proof.
\end{proof}
\begin{theorem}
    The following equation holds for all $a,b\in\mathbb{C}$
    \begin{align*}
        \vartheta_1(a+b)\vartheta_2(a-b)\vartheta_3\vartheta_4=\vartheta_1(a)\vartheta_2(a)\vartheta_3(b)\vartheta_4(b)+\vartheta_1(b)\vartheta_2(b)\vartheta_3(a)\vartheta_4(a).
    \end{align*}
\end{theorem}
\begin{proof}
    If $b$ is congruent to $0$ or $\frac{\pi}{2}$ modulo $(\pi,\pi\tau)$, then the theorem holds trivially by Lemma \ref{ThetaTransformLemma}. Therefore assume otherwise. Let us define the function
    \begin{align*}
        \psi_2(a)=\vartheta_1(a+b)\vartheta_2(a-b)\vartheta_3\vartheta_4
        -\vartheta_1(a)\vartheta_2(a)\vartheta_3(b)\vartheta_4(b)-\vartheta_1(b)\vartheta_2(b)\vartheta_3(a)\vartheta_4(a)
    \end{align*}
    By Lemma \ref{ThetaQPLemma}, we have $\psi_2(a+\pi)=\psi_2(a)$ and
    \begin{align*}
        \psi_2(a+\pi\tau)&=\vartheta_1(a+b+\pi\tau)\vartheta_2(a-b+\pi\tau)\vartheta_3\vartheta_4\\
        &-\vartheta_1(a+\pi\tau)\vartheta_2(a+\pi\tau)\vartheta_3(b)\vartheta_4(b)\\
        &-\vartheta_1(b)\vartheta_2(b)\vartheta_3(a+\pi\tau)\vartheta_4(a+\pi\tau)\\
        &=-q^{-2}e^{-4ia}\psi_2(a).
    \end{align*}
    Assume that $\psi_2$ is non-constant. Then by equation \eqref{GeneralN}, the number of zeros of $\psi_2(a)$ in a period parallelogram is $N=2$. We have $\psi_2(0)=\psi_2(-b)=0$ and furthermore
    \begin{align*}
        \psi_2\left( \frac{\pi}{2} \right)&=\vartheta_1\left( b+\frac{\pi}{2} \right)\vartheta_2\left( \frac{\pi}{2}-b \right)\vartheta_3\vartheta_4-\vartheta_1(b)\vartheta_2(b)\vartheta_3\left( \frac{\pi}{2} \right)\vartheta_4\left( \frac{\pi}{2} \right)\\
        &=0.
    \end{align*}
    But since $\psi_2$ only has two zeros in a period parallelogram, we must have that $b$ is congruent to either $0$ or $\frac{\pi}{2}$ modulo $(\pi,\pi\tau)$. Thus contradicting our assumptions. This completes the proof.
\end{proof}
Given complex numbers $a,b,c$ and $d$, we define the associated variables $a^*,b^*,c^*$ and $d^*$ by
\begin{gather*}
    2a^*=-a+b+c+d\\
2b^*=a-b+c+d\\
2c^*=a+b-c+d\\
2d^*=a+b+c-d.
\end{gather*}
Let us employ the convenience of notaion
\begin{align*}
    [j]=\vartheta_j(a)\vartheta_j(b)\vartheta_j(c)\vartheta_j(d) \qquad \text{and}\qquad [j]^*=\vartheta_j\left( a^* \right)\vartheta_j\left( b^* \right)\vartheta_j\left( c^* \right)\vartheta_j\left( d^* \right).
\end{align*}
We prove the first of the Jacobi's fundamental identities, from which the other $255$ may be deduced.

\begin{theorem}
    For every $a,b,c,d\in\mathbb{C}$, the following equation holds
    \begin{align*}
        2[3]=-[1]^*+[2]^*+[3]^*+[4]^*.
    \end{align*}
\end{theorem}
\begin{proof}
    Let us define the function
    \begin{align*}
        \psi_3(a)=2[3]+[1]^*-[2]^*-[3]^*-[4]^*.
    \end{align*}
    Under the translation $a\mapsto a+\pi$, we have the transformations
    \begin{gather*}
        [1]^*\mapsto -[2]^*,\qquad\qquad [2]^*\mapsto-[1]^*\\
        [3]^*\mapsto[4]^*,\qquad\qquad [3]^*\mapsto[4]^*
    \end{gather*}
    and therefore $\psi_3(a+\pi)=\psi_3(a)$. Similarly, under the translation $a\mapsto a+\pi\tau$, we have the transformations
    \begin{gather*}
        [1]^*\mapsto -q^{-1}e^{-2ia}[4]^*,\qquad\qquad [4]^*\mapsto -q^{-1}e^{-2ia}[1]^*\\
        [2]^*\mapsto q^{-1}e^{-2ia}[3]^*,\qquad\qquad [3]^*\mapsto q^{-1}e^{-2ia}[2]^*
    \end{gather*}
    and therefore $\psi_3(a+\pi\tau)=q^{-1}e^{-2ia}\psi_3(a)$. Assume that $\psi_3$ is non-constant. Then by equation \eqref{GeneralN}, the number of zeros of $\psi_3(a)$ in a period parallelogram is $N=1$. With a similar investigation of the roots of $\psi_3$, we obtain a contradiction. This completes the proof.
\end{proof}

\end{document}